\documentclass[a4paper,english, 11pt]{amsart}
\advance\oddsidemargin by -0.1cm
\advance\evensidemargin by -0.7cm
\advance\topmargin by -1.0cm 
\usepackage{amscd,amsmath,amsopn,amssymb,amsthm,multicol}
\usepackage[color,matrix, all, 2cell]{xy}
\usepackage{amscd}
\usepackage{verbatim}

 \usepackage{graphicx}
\usepackage[all]{xy}
\usepackage{slashed}
\usepackage{lscape}
 \usepackage{cancel}
\usepackage{mathrsfs}
\usepackage{euscript}
\usepackage{graphicx}
\usepackage{upgreek}
\usepackage{textgreek}
\usepackage{enumerate}
\usepackage{color}
\usepackage{tikz}
\usepackage{multirow}
\usepackage{cancel}
\usepackage{soul}
\usepackage{mathabx}
\usepackage{harmony}
\usepackage{comment}
\usepackage{tcolorbox}
\usepackage{marginnote}
\usepackage{wrapfig}
\usepackage{tikz}
 \numberwithin{equation}{section}

\DeclareMathOperator{\diag}{diag}
\DeclareMathOperator{\Ad}{Ad}
\DeclareMathOperator{\ad}{ad}

\DeclareMathOperator{\Lie}{\mathsf{Lie}}
\DeclareMathOperator{\Iso}{\mathsf{Iso}}

\DeclareMathOperator{\tr}{tr}

\DeclareMathOperator{\R}{\mathbb{R}}
\DeclareMathOperator{\Z}{\mathbb{Z}}

\DeclareMathOperator{\Cc}{\mathbb{C}}
\DeclareMathOperator{\Hh}{\mathbb{H}}

\DeclareMathOperator{\SO}{\mathsf{SO}}
\DeclareMathOperator{\Sp}{\mathsf{Sp}}
 
 \DeclareMathOperator{\Sl}{\mathsf{SL}}

\DeclareMathOperator{\G}{\mathsf{G}}

\DeclareMathOperator{\Ed}{\mathsf{End}}
\DeclareMathOperator{\GL}{\mathsf{GL}}

\definecolor{dblue}{rgb}{0.01,0.01,0.42}
\definecolor{red}{rgb}{0.57,0.11,0.15}

\newcommand{\fr}{\mathfrak}
\newcommand{\al}{\alpha}
\newcommand{\bee}{\beta}

\newcommand{\bb}{\mathbb}
\newcommand{\mc}{\mathcal}

\newcommand{\mf}{\mathsf}

\theoremstyle{plain}
\newtheorem{lemma}{Lemma} [section]
\newtheorem{theorem}[lemma]{Theorem}
\newtheorem{corol}[lemma] {Corollary}
\newtheorem{prop} [lemma]{Proposition}

\makeatletter
\newcommand{\tpitchfork}{%
  \vbox{
    \baselineskip\z@skip
    \lineskip-.52ex
    \lineskiplimit\maxdimen
    \m@th
    \ialign{##\crcr\hidewidth\smash{$-$}\hidewidth\crcr$\pitchfork$\crcr}
  }%
}
\makeatother

\theoremstyle{definition}
\newtheorem{definition}[lemma] {Definition}

\newtheorem{remark}[lemma] {Remark}
\newtheorem*{remark*}{Remark}

\def\bd{\begin{definition}}
\def\ed{\end{definition}}
\def\bt{\begin{theorem}}
\def\et{\end{theorem}}
\def\bl{\begin{lemma}}
\def\el{\end{lemma}}
\def\bp{\begin{prop}}
\def\ep{\end{prop}}
\def\br{\begin{remark}}
\def\er{\end{remark}}
\def\bc{\begin{corol}}
\def\ec{\end{corol}}
\def\be{\begin{equation}}
\def\ee{\end{equation}}

\def\rk{\operatorname{rk}}

\definecolor{dark}{rgb}{0.18,0.18,0.68}
\definecolor{mydark}{rgb}{1.08,0.08,0.08}
\definecolor{crew}{rgb}{0.2,0.5,0.2}
\definecolor{mmg}{rgb}{0.31,0.50,0.23}
\definecolor{dblue}{rgb}{0.01,0.01,0.44}
\definecolor{red}{rgb}{0.57,0.11,0.15}
\definecolor{cobalt}{rgb}{0.01, 0.31, 0.59}
\usepackage[colorlinks,citecolor=cobalt,linkcolor=cobalt,urlcolor=cobalt,pdfpagemode=UseNone,backref = page]{hyperref}

\usepackage{tcolorbox}
\definecolor{mycolor}{rgb}{0.122, 0.435, 0.898}

\begin{document}

\thispagestyle{empty}

\title{Reductive homogeneous Lorentzian manifolds} 

\author[Dmitri Alekseevsky, Ioannis Chrysikos and Anton Galaev]{Dmitri Alekseevsky, Ioannis Chrysikos and Anton Galaev}

\address{\hspace{-5mm}
{\normalfont Dmitri Alekseevsky, \ email: \ttfamily dalekseevsky@iitp.ru}\newline
Institute for Information Transmission Problems,
B. Karetny per. 19, 127051, Moscow, Russia and Faculty of Science,
University of Hradec Kr\'alov\'e,
Rokitanskeho 62, Hradec Kr\'alov\'e
50003, Czech Republic \newline\newline
{\normalfont Ioannis Chrysikos, \ email: \ttfamily ioannis.chrysikos@uni-hamburg.de}\newline
Department of Mathematics and Center for Mathematical Physics,
 University of Hamburg, 
Bundesstraße 55, D-20146 Hamburg, Germany\newline\newline
%
{\normalfont Anton Galaev, \ email: \ttfamily anton.galaev@uhk.cz}\newline
Faculty of Science,
University of Hradec Kr\'alov\'e,
Rokitanskeho 62, Hradec Kr\'alov\'e
50003, Czech Republic}

%
%
\begin{abstract}
We    study   homogeneous  Lorentzian manifolds  $M = G/L$  of  a  connected reductive Lie  group  $G$    modulo   a   connected reductive  subgroup  $L$, under   the assumption that   $M$   is (almost) $G$-effective and the isotropy representation is totally reducible.
We  show that  the  description of    such  manifolds   reduces to the case  of semisimple Lie groups $G$.
Moreover, we prove that  such a homogeneous space   is  reductive.
We  describe all  totally reducible  subgroups of the Lorentz  group and  divide them into three  types. The subgroups of \textsf{Type I}  are compact,  while the  subgroups of \textsf{Type II} and \textsf{Type III}  are non-compact. The  explicit   description  of the corresponding   homogeneous Lorentzian   spaces   of \textsf{Type II} and \textsf{III}  (under  some  mild  assumption) is given. We  also show that the   description  of Lorentz  homogeneous manifolds $M = G/L$ of  \textsf{Type I}    reduces  to the   description  of   subgroups $L$ such that $M=G/L$ is  an   admissible manifold, i.e.,  an  effective   homogeneous manifold that admits an invariant Lorentzian metric.  Whenever  the  subgroup  $L$  is  a  maximal    subgroup   with these properties, we call such a manifold  {\it minimal  admissible}.  We classify all   minimal   admissible homogeneous manifolds $G/L$ of  a compact  semisimple Lie  group $G$  and describe all invariant Lorentzian metrics on them.

\end{abstract}

\maketitle

\setcounter{tocdepth}{2}

\tableofcontents


\section*{Introduction}
This paper   is devoted to the investigation and classification   of  homogeneous $G$-effective Lorentzian manifolds
$M = G/L$ of  a reductive Lie group  $G$  modulo a reductive subgroup $L$.
Homogeneous  $G$-effective  Lorentzian manifolds   $M = G/L$ are naturally divided into two classes:  the class  of proper manifolds,  where  the action  of $G$ on $M$  is proper, or, equivalently, the  stabilizer  $L$ is  compact, and the class of non-proper homogeneous manifolds, where  $L$ is non-compact. 
A  remarkable result by N.~Kowalski \cite{Kow}  shows   that    the only  non-proper   homogeneous  Lorentzian manifold  $G/L$  of a simple  Lie group $G$    is    a space  of constant  curvature,  that   is  (up to a  covering)   the Minkowski space  or the (anti-)de\,Sitter  space.  
This result had been generalized by M.~Deffaf, K.~Melnick and A.~Zeghib \cite{DMZ} to the case of a  semisimple Lie group $G$: 
Any non-proper homogeneous Lorentzian manifold of a semisimple Lie group $G$ is a local
product of the  Minkowski  space  or the (anti-)de\,Sitter space and a Riemannian homogeneous manifold.

Our  starting point in this article is the classification  of    connected totally reducible  Lie   subgroups    of
the Lorentz group, which  we present in Section \ref{reductiveL}. There we show that totally reducible  subgroups  of the Lorentz group     fall into three types which we call \textsf{Type I, II} and \textsf{III}, respectively (see Theorem \ref{them1}).  \textsf{Type I}  groups  are compact  subgroups,  while \textsf{Type II} and \textsf{Type III}  consist of non-compact subgroups.
As a corollary, it is shown that   any homogeneous manifold $M=G/L$ of a reductive Lie group $G$  with   totally reducible isotropy subgroup admits a  reductive  decomposition. The  explicit   description  of (non-proper)    homogeneous Lorentzian  reductive  manifolds   of \textsf{Type II} and \textsf{Type III}  is given in Section~\ref{IIIII}, under some  mild  assumption.
Note that   all non-reductive (hence, non-proper) 4-dimensional homogeneous Lorenztian manifolds  had  been   classified  by Fels and Renner \cite{Fels}.
Proper  homogeneous Lorentzian manifolds $M = G/L$ of semisimple Lie groups $G$  were  studied   in \cite{Ale},  where     the structure  of  such  manifolds  had been elucidated. The   classification problem of such manifolds  is   equivalent  to the  description  of all admissible connected compact   subgroups $L$.  Recall that a connected Lie subgroup $L\subset G$ of $G$ is called {\it admissible}, if the homogeneous  manifold $M=G/L$ admits an invariant Lorentzian metric.  Since any  closed  subgroup  $L'\subset L$ is  also admissible, the
problem reduces to   the description  of  all {\it maximal admissible}  subgroups.  
Proper homogeneous Lorentzian manifolds $M=G/L$,  where  $L$ is a maximally admissible   subgroup  of a  semisimple Lie group  $G$  have  been  described in \cite{Ale}.  For this case  it  is  easy    to see   that a maximal admissible   subgroup  $L$   exists  only when $G$ is  simple.

In this paper   we   modify   the  notion  of   maximal admissible connected closed subgroups $L$  of  a Lie  group $G$, by  demanding   that  $M=G/L$  admits an invariant Lorentz  metric, such that there  is no  larger  closed connected   subgroup $\tilde{L} \supset L$  which makes the coset  $\tilde{M} = G/\tilde{L}$  an {\it effective} homogeneous Lorentzian manifold. 
Then,  the  manifold   $M = G/L$ is called  {\it minimal  admissible}. 
This  more  reasonable   definition of the notion of minimality, allows us to  extend the  results of \cite{Ale}.
In particular, in Section \ref{seccond}   we    give  a necessary  and sufficient condition   for  a  proper  homogeneous manifold  $M = G/L$ of  a  reductive Lie group $G$ to be a minimal  admissible  homogeneous space, in terms  of the   reductive decomposition
$\fr{g} = \fr{l} + \fr{m}$.
Then we   divide   such homogeneous spaces   into  three  subclasses \textsf{Ia}, \textsf{Ib}, and \textsf{Ic},
depending on whether  the   centralizer   $\fr{m}^{\fr{l}}$
is isomorphic  to $\R$ (\textsf{Ia}),  $\fr{m}^{\fr{l}}$ is three  dimensional simple regular Lie algebra $\fr{s}$  (\textsf{Ib}),   or  $\fr{m}^{\fr{l}}$ is the commutative Lie  algebra  $\R^k$  (\textsf{Ic}).
In Section  \ref{compactcase} we   use    the  above  result to describe  all  maximal   admissible  subgroups $L$ of a compact semisimple Lie group $G$   and classify all associated  minimal admissible   homogeneous manifolds $M= G/L$. Any  such  homogeneous space of \textsf{Type Ia}  is a   standard  homogeneous contact   manifold.
Any indecomposable minimal admissible  manifold  $M = G/L$ of  \textsf{Type Ib}   is the total space  of the  principal  $\Sp(1)$-bundle $$M = G/L \to  W = G/N_{G}(\fr{g}(\alpha))= G/\Sp(1)\cdot L $$ over  the    symmetric  quaternionic K\"ahler  manifold  $W$  corresponding to a simple  compact Lie group $G$ (a Wolf space), and hence  $M = G/L$  admits an invariant 3-Sasakian  structure. 
Finally we show that there are no  minimal  admissible    manifolds $M=G/L$   with compact $G$ of \textsf{Type Ic}.

In  Section  \ref{secnoncomp}  we generalize  the results of Section \ref{compactcase} to the case of non-compact   semisimple Lie groups $G$, and  describe   all proper  minimal  admissible  homogeneous   manifolds  $M = G/L$ with compact  stability   subgroup $L$ of\textsf{Type Ia}  and \textsf{Type Ib}.  In   particular, we show  that  such homogeneous space of \textsf{Type   Ia} (respectively \textsf{Type Ib})   are   non-compact  real  forms of   para 3-Sasakian manifolds (respectively,  pseudo 3-Sasakian manifolds). They  are   principal $\mathsf{SL}(2,\R)$-bundles   (respectively $\Sp(1)$-bundles)  associated  to the symmetric  para-quaternionic K\"ahler manifolds (or splittable quaternionic K\"ahler manifolds),  classified  in \cite{DJS} (respectively, the  symmetric pseudo-quaternionic K\"ahler manifolds, classified in \cite{AleC05}).

\section{Preliminaries}
In this article we adopt  the following notation and assumptions.
If the opposite is not stated, we    will    assume that  all Lie groups  are connected.
By a \textsf{homogeneous $G$-manifold}  $M = G/L$ we understand a  coset  space  of  a connected Lie group  $G$ modulo a closed connected   subgroup $L\subset G$ (called the  stability  subgroup). If the opposite  is not  stated,  we  will assume  that  the  homogeneous  manifold  $M = G/L$  is  almost effective.  This means  that       $G$ acts almost effectively  on $M = G/L$,   that is,   $L$ contains no non-discrete  normal subgroups of the Lie group $G$, or  equivalently, that the isotropy   representation $j: \fr{l} \to  \Ed(T_oM)$  of the  stability subalgebra
 $\fr{l}$ is exact.   
  Any  almost effective homogeneous $G$-manifold $M=G/L$  becomes  effective  if   we  factorize $G$  by the maximal (discrete) normal subgroup $\Gamma$ which acts  trivially  on  $M$.
 We  say that a Lie  subalgebra $\fr{l}$  of  a Lie   algebra $\fr{g}$  is  {\textsf{compact}} if the   adjoint  subalgebra  $\ad_{\fr{l}}|_{\fr{g}} \subset \ad_{\fr g}$  generates  a compact  subgroup $L$   of the  adjoint  group $\Ad_G$ corresponding to the Lie algebra  $\fr{g}$, and  it is called   \textsf{strongly  compact}  if  it generates a  compact   subgroup  of  the  simply connected  Lie group   $\tilde{G}$   with  $\Lie(\tilde{G})=\fr{g}$  (and hence,  of  any  Lie group ${G} $ with  $\Lie(G) =\fr{g}$).
 
\smallskip
Given a  homogeneous $G$-manifold $G/L$,  we will denote by $\fr{l} := {\textsf{Lie}}(L)\subset\fr{g} := {\textsf{Lie}}(G)$ the Lie algebras of $L\subset G$, respectively, and identify  the   tangent  space  $T_oG/L$ at the point  $o = eL$  with  the vector  space   $\fr{g}/\fr{l}$.  The   isotropy  representation   of the  stability group $L$ will be denoted by
 $j  : L \to  \GL(\fr{g}/\fr{l})$.
Recall that  a  homogeneous  space  $M=G/L$ is called \textsf{reductive}   if   the Lie algebra   $\fr{g}$  admits  a {\textsf{reductive decomposition}
\begin{equation} \label{Reductive decomposition} \fr{g}= \fr{l} + \fr{m},
\end{equation}
i.e., an $\Ad_L$-invariant decomposition  into a direct  sum  of  subspaces. Hence  we require that   the  subalgebra  $\fr{l}$  admits  an $\Ad_L$-invariant   complement $\fr{m}$ in $\fr{g}$,  which one can naturally identify with the tangent   space $T_oG/L$.  Then, the    isotropy  representation $j$ is identified  with the  adjoint  representation $\Ad_L |_{\fr{m}}$,    and the tangent bundle $TM$ is identified  with the associated vector bundle $G\times_{j}\fr{m}$.   It turns out that  $G$-invariant pseudo-Riemannian metrics bijectively correspond to
$\Ad_{L}$-invariant  pseudo-Euclidean  inner products on $\fr{m}$. In fact, since by assumption we work with connected stability groups, the $\Ad_{L}$-invariance can be equivalently replaced by $\ad_{\fr{l}}$-invariance.   Note that the  action  of a Lie group $G$ on   a    reductive homogeneous manifold  $M = G/L$ is   almost effective if and only if   the  isotropy representation   $\ad_{\fr{l}} |_{\fr{m}} \to   \fr{gl}(\fr{m})$  is exact,   and then  $\fr{l}$ has trivial intersection with the center $Z(\fr{g})$: $\fr{l}\cap Z(\fr{g}) = 0$. Next we  will  always  assume  this condition.

By a   \textsf{homogeneous  Lorentzian  $G$-manifold}  we understand a
  Lorentzian manifold  $(M,g)$   with a   transitive  almost  effective  isometric  action  of a  connected  Lie group   $G$. Then  we  identify  $M $
  with the   coset space   $M= G/L$,  where $L = G_o$ is the stability subgroup  of  a point $o \in M$. Below we  will always assume that  $G \subset \Iso(M,g)$ is a closed  subgroup  of the  full isometry group, and   that  $L$ is connected. Then   the isotropy  group $j(L)$ is a closed connected subgroup of the Lorentz  group.  Next we  are interested in the classification of    homogeneous  Lorentzian  manifolds $M = G/L$ with a reductive decomposition as in \eqref{Reductive decomposition} and with totally reducible isotropy representation.  To  simplify the exposition, we say  that  such  reductive decomposition   is \textsf{effective}.
   Then,  the  classification of   reductive  homogeneous Lorentzian manifolds with totally reducible isotropy representation   up to a covering,  reduces   to  the description
   of  all effective    reductive decompositions (\ref{Reductive decomposition}),  such that  the  isotropy  subalgebra
   $j(\fr{l}) := \ad_{\fr{l}}|_{\fr{m}}$  is   a totally reducible subalgebra  of   the Lorentz Lie algebra  $\fr{so}(\fr{m})$ (associated to some Lorentz  inner product   $g_{\fr{m}}$).}


 \section{Totally  reducible  subalgebras of the Lorentz  algebra}\label{reductiveL}
  In this section   we  classify  all  totally  reducible  subalgebras  of the Lorentz  algebra, and next specify the maximal ones.
In the following it is convenient to work with   the $(n+2)$-dimensional  Minkowski  vector space, which we shall denote by
\[
(V:=\R^{1, n+1}, \ g(u, v)=\langle u , v \rangle_{\mf{Min}}=-u_{0}v_{0}+\sum_{i=1}^{n+1}u_{i}v_{i})\,.
\]
The  corresponding connected {\textsf{Lorentz  group} will be denoted by $\SO^{0}(V)\equiv\SO^{0}(1, n+1)$ and we will write $\fr{so}(V)\equiv\fr{so}(1, n+1)$ for its Lie algebra, which is also referred to as the  \textsf{Lorentz  algebra}.  
Let us  fix a basis  $p, e_1, \ldots, e_n, q$  of $V$ (Witt basis), where $p, q$ are  isotropic vectors  with  $g(p,  q) =1$, and $e_1, \ldots, e_n$ is the orthonormal basis of the orthogonal complement to $\R p+\R q$ in $V$.  The  2-dimensional Minkowski subspace spanned by $p, q$ will be  denoted  by
	\[
	U=\R p+\R q,
	\]
	while   the vector space which is orthogonal   to $U$ will be denoted by $E$.  Obviously, $E$ is spanned by $e_1, \ldots, e_n$  and endowed with  the Euclidean metric $g|_{E}$, coincides with the Euclidean space~$\R^n$.
	\bd The direct sum decomposition
	\begin{equation}\label{eq1}
	V = \R p + E + \R q = U + E \,.
	\end{equation}
	is called a \textsf{standard decomposition} of $V$.
	\ed
	
	In terms of Lie algebras, we can  identify  skew-symmetric  endomorphisms  with bivectors, and so the Lorentz   algebra $\fr{so}(1, n+1)\equiv\fr{so}(V)$ with $\Lambda^{2}V$.
	Then,
	\bl The Lorentz algebra   $\fr{so}(V)$ admits a {\sf depth-1 grading}
	\[
	\fr{so}(V) = \fr{g}^{-1} + \fr{g}^0 + \fr{g}^1 = q\wedge E  + \big(\R p\wedge q + \fr{so}(E)\big) + p \wedge E\,,
	\]
	which is  the  eigenspace decomposition  of the endomorphism $d = \ad_{p \wedge q}$.
	\el
	\begin{proof}
		This claim immediately follows by  (\ref{eq1}) and the identification  $\fr{so}(V)=\Lambda^{2}V$.
	\end{proof}
Note that 	the   parabolic subalgebra
	\[
	\fr{p}:=\fr{g}^{0}+\fr{g}^{1}= \R p \wedge q + \fr{so}(E)  + p \wedge E\,,
	\]
	coincides with the Lie  algebra  of the Lie group ${\SO^0(V)}_{\mathbb{R}p}$ preserving the isotropic line $\mathbb{R}p$.
	
	Consider now an orthogonal direct sum decomposition
	\[
	E = H + H^{\perp}
	\]
	of the Euclidean vector space  $(E\cong\R^n, g|_{E})$ such that  $\dim H = k, \dim H^{\perp} = n-k$, $0\leq k\leq n-1$  and let us denote by
	\begin{equation}\label{vH}
	V(H):=\R p+H+\R q\subseteq V
	\end{equation}
	the  Lorentzian vector space  of dimension $k+2$.  Note that $V$ can be expressed as the direct sum decomposition
	\begin{equation}\label{HvH}
	V=(\R p+H+\R q)+H^{\perp}=V(H)+H^{\perp}\,.
	\end{equation}
	Moreover,
	\begin{itemize}
		\item for $k=n$ we have $V(H)=V$, while
		\item for $k=0$, $E=H^{\perp}$ and $V(H)=U=\R p+\R q$.
	\end{itemize}
	Set    $\fr{h}_{k}:= \fr{so}(V(H)) + \fr{so}(H^{\perp})$. Then $\fr{h}_{k}$  is    a direct sum of a Lorentz  algebra and an orthogonal Lie algebra, in particular
	\begin{equation}\label{hk}
	\fr{h}_{k}=\fr{so}(V(H)) + \fr{so}(H^{\perp})\cong \fr{so}(1, k+1)+\fr{so}(n-k).
	\end{equation}
	
	\bp  \label{maxsuba} A  maximal subalgebra $\fr{h}$ of $\fr{so}(V)$ is conjugated  to one of the  following subalgebras:
	\begin{itemize}
		\item the  maximal compact subalgebra $\fr{so}(n+1)$;
		\item $\fr{h}_{k}=\fr{so}(V(H)) + \fr{so}(H^{\perp})\cong \fr{so}(1, k+1)+\fr{so}(n-k)$, $1 \leq k\leq n-1$;
		\item the parabolic  subalgebra $\fr{p} = \R p \wedge q + p \wedge E + \fr{so}(E)$
	\end{itemize}
	\ep
	\begin{proof}
		It is clear that the algebras listed above are maximal subalgebras of $\fr{so}(V)$.
		Let $\fr{h}$ be a maximal subalgebra of $\fr{so}(V)$. It is well-known that there are no proper irreducible subalgebras in $\fr{so}(V)$. Hence $\fr{h}$ preserves a proper subspace of $V$.
		If
		$\fr{h}$ preserves a proper non-degenerate subspace of $V$, then $\fr{h}$ preserves also a Minkowski subspace $W\subset V$. If $\dim W=1$, then $\fr{h}$ is conjugated to   $\fr{so}(n+1)$. 
		If $\dim W=2$, then $\fr{h}$ is conjugated to the parabolic subalgebra $\fr{p}$ (note that $\fr{so}(1,1)+\fr{so}(n)$ is not a maximal subalgebra of $\fr{so}(1,n+1)$). Finally,  if $\dim W\geq 3$, then $\fr{h}$ is conjugated to $\fr{h}_{k}$ for some  $1 \leq k\leq n-1$. 
			Suppose now that $\fr{h}$ does not preserve any non-degenerate subspace of $V$. Then $\fr{h}$ preserves a degenerate subspace $W\subset V$ and  the isotropic line $W\cap W^\bot$. Let $p$ be a generator of $W\cap W^\bot$ and fix  an isotropic vector $q\in V$ as above. Then we obtain the decomposition
			$V=\mathbb{R}p+E+\mathbb{R}q$,
			 and conclude that
		\[
		\fr{h}=\mathbb{R}p\wedge q+\fr{so}(E)+p\wedge E=\fr{p}\,.
		\]
		This completes the proof.
			\end{proof}
	
	Let us now pose the main theorem of this section.
	\bt \label{them1} A  totally reducible subalgebra  $\fr{h}$ of   the  Lorentz  algebra  $\fr{so}(V)$
is  conjugated  to a subalgebra of one of the following types: \begin{itemize}
\item[\textsf{Type  I}:]  a subalgebra  of   the maximal  compact Lie  algebra   $\fr{so}(n+1)$;
\item[\textsf{Type II}:] a  subalgebra  of  the  form   $\fr{so}(W) + \fr{k}$, where $W=V(H)$ is  a Lorentzian  subspace   of  $V$, $\dim W\geq 3$, and   $\fr{k}$ is  a  compact  subalgebra  of  $\fr{so}(W^{\perp})$;
\item[\textsf{Type III}:]   a  subalgebra  of  the  form  $\R d   \oplus \fr{k}$, 
where   $d = p \wedge q  + C_0$ with $C_0 \in \fr{so}(E)=\fr{so}(n)$ (which is possibly zero),   and  $\fr{k}$ is  a  subalgebra  of  the centralizer $C_{\fr{so}(n)}(C_0)$.
\end{itemize}
\et

\begin{proof} Note that the Lie algebra $\fr{so}(V)$ is of \textsf{Type II}.
	First suppose that $\fr{h}$ is different from $\fr{so}(V)$ and $\fr{h}$ does not preserve any proper non-degenerate subspace of $V$. By Proposition \ref{maxsuba}, $\fr{h}$ is conjugated to a subalgebra of the parabolic algebra and it preserves the line $\R p$. Since $\fr{h}$ is totally reducible, it preserves a  vector subspace $V_0\subset V$ complementary to $\R p$. By the assumption, $V_0$ is degenerate. Consequently $\fr{h}$ preserves the isotropic line $\ell=V_0\cap V_0^\bot$, and the vector subspace $\R p+\ell\subset V$. It is clear that the vector subspace $\R p+\ell\subset V$ is Lorentzian, and this gives a contradiction. 
	
\noindent	Suppose now that $\fr{h}$ preserves a non-degenerate subspace of $V$. Then  $\fr{h}$ preserves a Lorentzian subspace $W\subset V$, and we may assume that $\fr{h}$ does not preserve any proper non-degenerate subspace of $W$. As we have seen just above, this implies that the projection of $\fr{h}$ to $\fr{so}(W)$ coincides with  $\fr{so}(W)$.
	If $\dim W=1$, then $\fr{h}$ is of \textsf{Type I}. If If $\dim W=2$, then $\fr{h}$ is of \textsf{Type III}. If $\dim W\geq 3$, then $\fr{h}$ is of \textsf{Type II}.
	\end{proof}

\bc Connected totally reducible subgroups of the Lorentz group $\SO(V)$ are divided into the following 3 types:
\begin{itemize}
	\item[\textsf{Type  I}:]  connected subgroups  of   the maximal  compact Lie  group   $\SO(n+1)$;
	\item[\textsf{Type II}:] connected subgroups of the form   $\SO(W) \cdot K$, where $W=V(H)$ is  a Lorentzian  subspace   of  $V$, $\dim W\geq 3$, and   $K$ is  a connected subgroup  of  $\SO(W^{\perp})$;
	\item[\textsf{Type III}:]    subgroups  of  the  form  $T\cdot K$,  
	where   $K\subset\SO(E)$ is a  connected subgroup and $T$ is a 1-parameter subgroup of $\SO(1,1)\times\SO(E)$ not contained in $\SO(E)$ and commuting  with~$K$.
\end{itemize}
\ec

\begin{remark}
	It is not hard to show that reductive subalgebras $\fr{h}$ of $\fr{so}(V)$ are exhausted by the subalgebras specified in   Theorem \ref{them1}, together  with subalgebras $\fr{h}$ of the following form.
	Consider an orthogonal decomposition $E=E'+E''$ and the corresponding decomposition $$V=(\R p+E'+\R q)+E''.$$
	Let $\fr{k}\subset\fr{so}(E'')$ be a subalgebra and let $\varphi:E'\to\fr{so}(E'')$ be a linear map such that $\varphi(E')\subset\fr{so}(E'')$ is commutative, it commutes with $\fr{k}$ and has trivial intersection with $\fr{k}$. Then, the Lie subalgebra $$\fr{h}=\{p\wedge X+\varphi(X)|X\in E'\}+\fr{k}\subset\fr{so}(V)$$
	is reductive, but not totally reducible.	
	\end{remark}

	As it is already indicated, the above theorem motivates  us to introduce the a  following definition.
	\bd
	$(\al)$ We refer   to  the   subalgebras specified by Theorem \ref{them1}    as the  \textsf{totally reducible  Lie   subalgebras of the Lorentz algebra $\fr{so}(V)=\fr{so}(1, n+1)$} of  \textsf{Type I,  II}, and \textsf{III}, respectively,  and similarly for the corresponding connected closed Lie subgroups of $\SO^{0}(V)\equiv\SO^{0}(1, n+1)$. \\
	$(\bee)$ A homogeneous  Lorentzian manifold $M = G/L$ is called a $G$-manifold   of \textsf{Type I} (respectively, \textsf{II}, \textsf{III}) if   the   isotropy  Lie algebra $j(\fr{l})$  is of   \textsf{Type I} (respectively, \textsf{II}, \textsf{III}).
	\ed
	Note that      an effective   homogeneous Lorentz  manifold $M = G/L$ of   \textsf{Type I} (respectively, \textsf{II}, \textsf{III})  has a  compact (respectively, non-compact) stabilizer  $L$.	Let us now prove the following result.

\bp  
Each homogeneous Lorentzian manifold  $M = G/L$  of a reductive Lie group  $G$ with a totally reducible isotropy subgroup   is   a reductive   homogeneous manifold.
 \ep

\begin{proof} We  may assume that  the   homogeneous  Lorentz manifold $M = G/L$ is effective.  Then the  isotropy representation $j$ is exact.   Since    the  isotropy  group $j(L)$  of  \textsf{Type I}  is   compact  and the isotropy group  of \textsf{Type II}  is an almost direct product $ j(L) = \SO^0(V(H)) \cdot K$  of a (simple)  Lorentz Lie group   and a  compact Lie group $K$,    the  adjoint   representation  $\Ad_L |_{\fr{g}}$  is totally reducible. Moreover, there is an     $\Ad_L $-invariant subspace $\fr{m}$ complementary  to $\fr{l}$.  It  remains  to  check that   the isotropy  group  $j(L) = T \cdot K  $, geberated  by the   isotropy Lie  algebra $j(\fr{l}) = \R d \oplus \fr{k}$ admits a reductive decomposition.   Since $\fr{k} \subset \fr{g}$  is a compact  subalgebra, we may  chose  an $ j(\fr{k})$-invariant      subspace
 $\fr{m}$  complementary   to  $\fr{l}$   and  write   the  Lie   algebra $\fr{g}$ in the form
 $$ \fr{g} = \fr{l} + \fr{m} = \fr{l} + \R p + \R q + E_0 + E'$$
 such that $j(d) = p \wedge q + C_0$,
 $$C_{\fr{g}}(\fr{k}) = \R d +Z(\fr{k}) + \R p + \R q + E_0,  $$
   and $j(\fr{k})E' = E'$.
 Since  $j(d)p, \, j(d)q  \in C_{\fr{g}}(\fr{k}) $,   we  have
 $$ j(d) p = p + \alpha d + z_1, \quad  j(d) q = -q +\beta d + z_2, \quad z_1, z_2  \in Z(\fr{k}). $$
 Taking $p': = p +\alpha d + z_1$, $q' := q - \beta d - z_2   \in C_{\fr{g}}(\fr{k})$, we get  the $j(\fr{l})$-invariant 
subspace  $\fr{m}' := \R p' + \R q' + E_0 + E'$ 
   complementary   to $j(\fr{l})$ and such that 
$j(d)p' = p',$   $j(d)q' = -q'$.
 \end{proof}

\section{Admissible homogeneous manifolds}\label{AdmissibleHomMan}
In this section we   introduce  a small  but important  modification  of the
 the notion of  \textsf{admissible subgroups}, introduced by the first  author in \cite{Ale}.
\bd\label{def1}
 Let $G$ be a connected Lie group. Then
\begin{itemize}
	\item[\textsf{(i)}]  A closed connected  subgroup $L$  of    $G$ is  called \textsf{admissible}
 if  the  homogeneous  space  $M = G/L$ is   almost  effective   and admits  an invariant Lorentz  metric~$g$. Then  $M=G/L$ is said to be an   \textsf{admissible homogeneous space}.
 
  \item[\textsf{(ii)}] An admissible  subgroup $L\subset G$ is said to be  \textsf{maximal  admissible}    and  then the manifold $M = G/L$  is said to be   a \textsf{minimal  admissible}, if  there  is no  larger  almost effective admissible  subgroup $\tilde{L}\supset L$ with  $\dim  \tilde{L} > \dim L$.
  
  \item[\textsf{(iii)}] The Lie  algebra $\fr{l}=\Lie(L)$  of a  (maximal)  admissible   subgroup $L \subset G$  is called  a (maximal) admissible  subalgebra  of $\fr{g}$.
  \end{itemize}
\ed

\br  Definition \ref{def1} provides  a  modification  of the  notion of a minimal    admissible  homogeneous manifold  given in \cite{Ale}, where   the   effectivity condition of  the  action  of $G$  on   $G/\tilde{L}$  is not required.
This  modification  is  useful  by the  following reason.
Let  $M_1 = G/L_1$  be an almost effective minimal  admissible  manifold. Consider the direct product of $M_1$ with the sphere,  $M = G/L= G_1/L_1\,\times\SO(m+1)/\SO(m)$. Then $M$  is a minimal   admissible  manifold. However, it is not a minimal  admissible manifold  in the  sense of \cite{Ale}.
  Indeed,  $L_1 \times\SO(m+1)$ is an admissible  (but not   locally effective)  subgroup of $G = L_1 \times\SO(m+1)$, which contains  $L = L_1 \times\SO(m)$.
   Below  we  will construct    a  large new series  of   minimal  admissible homogeneous  manifolds  $G/L$    of  compact semisimple Lie groups $G$.
\er

  \subsection{Admissible homogeneous manifolds   of reductive  non-semisimple  Lie groups}

 Let  $M = G/L$  be  a  homogeneous   manifold  of  a  reductive Lie  group  $G = G^s \cdot Z(G)$,   where  $G^s$ is the  semisimple  part  of $G$,   $Z(G)$ is  the   connected non-trivial center of $G$,  and  $L\subset G$ is a reductive  subgroup. Then   the Lie  algebra
   $\fr{g} = \fr{g}^s + Z(\fr{g})$   is a direct  sum of the    semisimple  ideal $\fr{g}^s$ and  the center $Z(\fr{g})$ of $\fr{g}$.  Due   to the effectivity it holds $\fr{l} \cap Z(\fr{g}) =0$.
  Denote by
\[
{\rm pr}_{\fr{g}^s} : \fr{l} \to  \fr{l}_s := {\rm pr}_{\fr{g}^s}(\fr{l}) \subset \fr{g}^s
\]
the natural projection  of   $\fr{l}$ to  $\fr{g}^s$, which is an isomorphism   between the Lie  algebras $\fr{l}$ and $\fr{l}_s$.  The  adjoint   actions  of $\fr{l}$ and $\fr{l}_s$ on $\fr{g}$ coincide.
             Since   $\fr{l}_s$   is  reductive,  there   is  a  reductive  decomposition $\fr{g}^s = \fr{l}_s    + \fr{m}_{s}$   of $\fr{g}^s$. This gives the reductive decomposition 
 \[
   \fr{g} = \fr{l} + \fr{m} =\fr{l} + (\fr{m}_{s} + Z(\fr{g}))\,.
   \]
               The centralizer $\fr{m}^\fr{l}$ is given by
\[
 \fr{m}^{\fr{l}} = \fr{m}_s^{\fr{l}} + Z(\fr{g}) \neq 0\,.
 \]
             As a consequence we obtain the following.
           
           \bt Let   $M = G/L$ be  a  homogeneous    manifold   of a reductive group  $G = G^s \cdot Z(G)$  with a  non-trivial connected center $Z(G)$  and  a reductive stability subgroup $L$,  and let  $M^s = G^s/L_s$ be  the  associated    homogeneous manifold  of  the semisimple  Lie   group $G^s$ described above.             Then
             any invariant  Riemannian  or Lorentzian metric   on  $M^s = G^s/L_s$ may be  extended to an invariant  Lorentz metric  on  $M = G/L$.  \et
             
           \begin{proof}
            Let  $\fr{g} = \fr{l} + (\fr{m}_s   + \fr{z})$ be  a  reductive decomposition as above.  Then  any $j(L)$-invariant     Euclidean  or Lorentzian  inner product  $g_{\fr{m}_s}$ in   $\fr{m}_s$   is  extended  to the  $j(L)$-invariant Lorentzian   inner product $g_{\fr{m}} = g_{\fr{m}_s} \oplus g_{Z(\fr{g})}$ in $\fr{m}$,  where  $g_Z(\fr{g})$   is  any Lorentzian or Euclidean   inner product  in~$Z(\fr{g})$.
           \end{proof}

  \subsection{Proper admissible  homogeneous manifolds  }
  Let $M=G/L$  be  a  proper  homogeneous manifold, i.e., the stabilizer  $L$ is  compact.      We   fix a reductive decomposition
  and  identify  the isotropy   representation $j : L \to\GL(T_oM) = \GL(\fr{m})$   with the    restriction
    $\Ad_L|_{\fr{m}}$  to $\fr{m}$ of the  adjoint  representation   $\Ad_L|_{\fr{g}}$.
  Since  the  compact   linear group $\Ad_L \subset\GL(\fr{g})$ is totally  reducible,   there  exists   an $\Ad_L$ invariant   complement $\fr{m}$  to $\fr{l}$ in $\fr{g}$.   The homogeneous space $M$   admits  an  invariant Riemannian metric $g_M$,  which   is     determined  by  an $j(L)$-invariant  Euclidean metric   $g_{\fr{m}}$  on $\fr{m}$.
 Such a metric  $g_{\fr{m}}$  can be  defined  as  the barycentre  of   the compact   orbit   $j(L)^*g_0 \subset \rm{Met}(\fr{m})$ of the  natural action  of the group $j(L)$ in  the   convex  cone   $\rm{Met}(\fr{m})$ of  Euclidean  inner products  on $\fr{m}$, where  $g_0$ is a fixed Euclidean  metric.      Assume  that $j(L)$   preserves  a  line  $\R Z \subset \fr{m} $ on $\fr{m}$.   Then,    for  sufficiently  big  number $\lambda$,   
 \[
 g_{\lambda}:= g_{\fr{m}}  - \lambda  Z^* \otimes Z^*\,,\quad \quad   Z^* := g_{\fr{m}}\circ Z
 \]
    defines  an $j(L)$-invariant  Lorentzian  scalar product,  which  defines  an  invariant Lorentz metric on $M = G/L$. Conversely, any  invariant Lorentzian  metric  can be obtained  by this  construction. 

\bp\label{staolderes}\textnormal{(\cite{Ale})}
Let $M=G/L$    be a  homogeneous manifold  with compact connected stability  subgroup  $L$.
 The  manifold $M = G/L$  is admissible  if and  only if   the  centralizer
 $\fr{m}^{\fr{l}}$ is non-trivial.
  Moreover, let  $g_{\fr{m}}$   be   an $j(L)$-invariant  Euclidean metric  in $\fr{m}$,  and  $0 \neq Z \in \fr{m}^{\fr{l}}$.   Then  for   sufficiently  large  number  $\lambda$,  the    formula $g_{\lambda}:=  g_{\fr{m}}   - \lambda Z^* \otimes Z^*$, where  $Z^* := g_{\fr{m}}\circ Z$,
   defines   a Lorentzian $j(L)$-invariant   inner  product  on $\fr{m}$, which is extended  to a   $G$-invariant Lorentzian metric    on  $M$.
\ep

\bc   If $M= G/L$  is  an  admissible  homogeneous manifold  with compact stabilizer,  then  any closed subgroup  $L' \subset L$  is   admissible, that is, the  manifold $M'= G/L'$
admits  a $G$-invariant Lorentzian metric.
\ec
\begin{proof}
 Denote  by $\fr{l} = \fr{l}' + \fr{q}$ a reductive decomposition  for   $L/L'$.
 Then  $\fr{g} = \fr{l}' + \fr{m}' = \fr{l}' + (\fr{q} + \fr{m})$
is  a reductive decomposition for   $G/L'$ and the  action of  $G$   on $G/L'$ is almost effective. It holds $C_{\fr{m}'}(\fr{l}')\supset C_{\fr{m}}(\fr{l})$. Then, by Proposition \ref{staolderes} it follows that $M'$ is  an  admissible  homogeneous  manifold of $G$.
\end{proof}

\subsection{Condition for  minimal  admissibility}\label{seccond}
Next our aim is to describe sufficient   conditions  for  a    homogeneous manifold $M = G/L$ of a reductive group $G$   with compact stability subgroup $L$     to be a minimal   admissible homogeneous space.

\bl\label{key lemma} (Key lemma)  Let  $M = G/L$   be  a minimal  admissible  homogeneous    manifold  of   a reductive  Lie group $G$   with a compact stabilizer   $L$. 
Then, the   reductive  decomposition  can   be chosen in such a way that  
\[
 \fr{g}  =\fr{l} + \fr{m}= \fr{l} + \fr{m}^{\fr{l}} + \fr{m}', \quad
\fr{m}' = [\fr{l}, \fr{m}]
\]
 and $\fr{m}^{\fr{l}}$   is   a  \textsf{ reductive  subalgebra} of $\fr{g}$.
 Moreover, there  are  three  possibilities: \begin{itemize}
 \item[$(a)$] $\fr{m}^{\fr{l}} = \R Z$   is a compact  1-dimensional Lie  algebra which is not  in the center of $\fr{g}$,
 and  $ C_{\fr{g}}(Z) = \fr{l} + \R Z$.
  \item[$(b)$]  $\fr{m}^{\fr{l}} = \fr{s}$ is  a rank-one   regular   simple Lie  algebra, i.e., $\fr{sp}(1)$   or  $\fr{sl}(2,\R)$,  such that $C_{\fr{g}}(\fr{s})=\fr{l}$.
     For   any element $Z \in \fr{s}$  which defines  a compact subalgebra  $\R Z$,    it holds $C_{\fr{g}}(Z) =\fr{l} + \R Z$.
  \item[$(c)$]  $\fr{m}^{\fr{l}}= \R^k$ is a commutative  subalgebra    which does not contain any  compact non-central  subalgebra  $\R Z$.
\end{itemize}
  For any  subalgebra of type  (a) and (b),  the  normalizer
   $N_{\fr{g}}(\fr{m}^{\fr{l}})  = \fr{l} \oplus  \fr{m}^{\fr{l}}$
   of the   subalgebra  $\fr{m}^{\fr{l}}$  is a   reductive    subalgebra of $\fr{g}$ of maximal rank.
     \el

 \begin{proof}
 	 	Let
 	$   \fr{g} = \fr{l} + \fr{m}$
 	be a  reductive  decomposition associated  to  the   manifold~$M$.
 	It is clear that for the centralizer $C_{\fr{g}}(\fr{l})$ 
 	of   the   stability  subalgebra $\fr{l}$ in   $\fr{g}$ it holds 	$ C_{\fr{g}}(\fr{l})  = Z(\fr{l}) + \fr{m}^{\fr{l}}$,
 	where  $\fr{m}^{\fr{l}}= C_{\fr{m}}(\fr{l})$ is non-trivial by Proposition~\ref{staolderes}. 
 	Since the centralizer  of  a   compact Lie  algebra in  a  reductive Lie
algebra  is  reductive,  $C_{\fr{g}}(\fr{m}^{\fr{l}})$  is a  reductive Lie algebra. Since $Z(\fr{l})\subset C_{\fr{g}}(\fr{m}^{\fr{l}})$ is an ideal, there is  a complementary reductive  ideal, which we denote again by $\fr{m}^{\fr{l}}$.   We   may  extend    $\fr{m}^{\fr{l}}$   to an $\ad_{\fr{l}}$-invariant     subspace $\fr{m}$ complementary  to $\fr{l}$.  Then, the  new reductive  decomposition  $\fr{g}=\fr{l} + \fr{m}$  satisfies the  desired properties.

  \noindent Now we prove that if   the  reductive  Lie  algebra $\fr{m}^{\fr{l}}$    has  a compact  non-central  subalgebra   $\R Z$,   then  it  has rank  one.
   Indeed, the Lie subalgebra $\tilde{\fr{l}}:= \fr{l} +\R Z$ is   compact and  effective  if  $Z$ is not   in the center  of $\fr{g}$.  If  $\mathrm{rk} (\fr{m}^{\fr{l}}) >1$,   then  there  is an element
   $Z' \in   \fr{m}^{\fr{l}}$  non-proportional to $Z$, which  commutes  with $Z$, and hence
   with $\tilde{\fr{l}}$.     
   Then, $\tilde{\fr{l}}$ generates a   compact   admissible  subgroup of $G$, a fact which   contradicts  the maximality of $L$. 
Therefore,  if  the  reductive Lie algebra $\fr{m}^{\fr{l}}$ has a  compact non-central  subalgebra  $\R Z$, then  it  is  either  isomorphic to $\R Z$ or to a simple rank-one  subalgebra, and hence to $\fr{sp}(1)$
or to $\fr{sl}(2, \R)$.  If  there  is no compact non-central  subalgebra $\R Z$, then the   reductive  subalgebra $\fr{m}^{\fr{l}}$  is commutative.

 \noindent  Let us now  check that  if   $\R Z \subset \fr{m}^{\fr{l}}$  is   a compact non-central subalgebra, then the normalizer satisfies $N_{\fr{g}}(Z) = \fr{l} + \R Z$. The normalizer can  be   decomposed   as
\[
 N_{\fr{g}}(Z) = \fr{l} +\R Z + C_{\fr{m}'}(Z)\,,
 \]
    where  $\fr{m}' := [\fr{l}, \fr{m} ]$ is the   complementary  to  $\fr{m}^{\fr{l}}$ subspace of $\fr{m}$.  If $ C_{\fr{m}'}(Z)\neq 0 $, then  the compact  subgroup $\tilde{L}$ generated by  $\tilde{\fr{l}} = \fr{l} + \R Z$   will be    admissible, which is impossible.\end{proof}

\section{ Minimal compact    homogeneous  Lorentz manifolds   with  compact  stabilizer}\label{compactcase}
  Here   we describe   compact    homogeneous Lorentzian manifolds  $M = G/L$   with compact stabilizer. Note that   for   a compact  group  $G$    there  is no minimal admissible manifolds $M=G/L$ of  \textsf{Type  Ic}.  Hence, below we will   discuss  only  the  homogeneous manifolds  of  \textsf{Type Ia}   and \textsf{Ib},
    and  extend   the results of \cite{Ale}.

\subsection{\textsf{Type Ia}}
The  classification  problem for the     \textsf{Type Ia}   homogeneous manifolds  $M =G/L$  of a  compact  semisimple  group $G$ reduces  to  the description  of all
 effective reductive  decompositions  of  a compact  semisimple Lie  algebra $\fr{g}$ of the form
\begin{equation} \label{standard decomposition}
  \fr{g }= \fr{l} + \fr{m}= \fr{l} + (\R Z   +\fr{m}')\,,
  \end{equation}
 where $ \fr{h} = \fr{l} + \R Z $ (direct sum)  is the  centralizer   of  a non-trivial  element $Z \in \fr{m}$, that is, $\fr{h}=C_{\fr{g}}(Z)$,  such that  the   subalgebra $\R Z$ is compact. 
 Without loss  of generality, we may assume   that $Z$ is $B$-orthogonal  to $\fr{l}$  and $B(Z,Z) =-1$,  where  we denote by $B$   the Killing  form.
 In    \cite{AS} the   vector $Z \in \fr{g}$   is called the {\it contact  element},  and it  is shown that it generates a  closed  subgroup  $T^1_Z$.

   Let  $H = C_{G}(T^1_Z) $ be  the  (closed)  subgroup of $G$, generated  by  $\fr{h}$. It  is  the  centralizer  of the    1-parameter  subgroup   $T^1_Z = \overline{\exp \R Z}$.  Then  the  homogeneous space $F = G/H$ is a so-called {\it (generalized) flag manifold}. It is well-known that in this case  $H$ is connected. Since  $\fr{h} = \fr{l} +  \R Z$, the group  $H = L\cdot T_Z^1  $ is a  product  of  the normal subgroups  $L, \, T^1_Z$  and the  intersection $L \cap T^1_Z$  is  a  finite cyclic  normal central   subgroup. 
   Hence the  natural projection  
   \[
   \pi: M = G/L \to F = G/H
   \]
     is a circle  bundle  over   the flag manifold $F $. 
  Moreover, the $ \Ad_{L}$-invariant  subspace $\fr{m}'$  defines   the  invariant contact  structure $\mc{D} \subset TM$ on   $M = G/L$.   The corresponding invariant   contact 1-form $\theta $   is the  extension  of the  1-form $\theta_o:= B \circ Z$.

   \bd \cite{AS}  The manifold $M = G/L$ which is  the total space of the   principal circle bundle
   $\pi : M = G/L \to F = G/H$ over the  flag manifold $F = G/H$, endowed   with the  invariant contact  structure   $\mc{D}\subset TM$   associated to   the subspace  $\fr{m}'$ of {\rm (\ref{standard decomposition})},   is called  a  {\it standard  homogeneous contact manifold}. Moreover,  $\mc{D}$
  is referred to as   the {\it standard homogeneous contact structure}.
  \ed

  For  a given  compact  semisimple Lie group $G$,    standard homogeneous contact  manifolds $M = G/L$   are  in  one-to-one  correspondence  with  elements $Z \in \fr{g}$ (defined  up to  a scaling),      generating     closed  1-parameter  subgroups $T^1_Z = \overline{\exp \R Z}$ of $G$.  
 As above,  denote   by $\fr{h} = C_{\fr{g}}(Z)  = \fr{l} + \R Z$  the   centralizer of $Z$ in $\fr{g}$, where $\fr{l}$ is the $B$-orthogonal complement   of $\R Z$ in $\fr{h}$. In addition,  consider  the $B$-orthogonal standard  decomposition \eqref{standard decomposition} associated  to $Z$. In such terms one can prove   the following

     \bt \cite{AS} The subalgebras $\fr{l}$ and $\fr{h}= \fr{l} +  \R Z$ (direct sum)  generate   closed  subgroups $L$ and $H$ of $G$,
   and  the  natural $G$-equivariant  projection $ \pi : M= G/L \to F = G/H$
   is   a principal $  T^1_Z$-bundle   over  the   flag manifold $F = G/H$, 
    with the invariant contact  structure associated to the  subspace  $\fr{m}'$.
    The   action  of $G$ on $M = G/L$ is effective if  $G$ has trivial center. 
     \et
     Now  we are ready to present   the  theorem  which  describes  all compact homogeneous Lorentzian manifolds  $M = G/L$  of    \textsf{Type  Ia}.
     
     \bt
      Let  $G$ be a semisimple  compact  Lie group.  
        Then   all minimal  admissible  homogeneous  $G$-manifolds  are exhausted by the   standard homogeneous  contact  manifolds  $M = G/L$, which are  the circle  bundles  $\pi : G/L \to  F = G/H$ over    flag manifolds  associated to the   standard   decomposition~\eqref{standard decomposition}.   Let $g_{\fr{m}'}$ be  any
      $\ad_{\fr{l}}$-invariant Euclidean  metric  on $\fr{m}'$. Then any  invariant Lorentzian metric  on $M = G/L$  is   the invariant  extension of the  pseudo-Euclidean  metric  on  $T_oM= \fr{m} = \R Z + {\fr{m}}'$ of the  form
      $g_{\fr{m}}  = -\lambda Z^* \otimes Z^* +g_{\fr{m}'}$, where  $Z^*:= B \circ Z$ and  $\R\ni\lambda >0$. 
     \et
   \begin{proof}  It is sufficient  to show that  $M= G/L$ is  a minimal  admissible homogeneous  manifold. It follows  from the  relation $C_{\fr{g}}(\fr{l}) = \fr{l} + \R Z$, which implies  $[\fr{l}, \fr{m}']= \fr{m}'$  and
   $C_{\fr{g}}(Z) = \fr{l} + \R Z$. Using these relations,  one can check   that  any admissible  Lie   algebra
\[
\tilde{\fr{l}} = \tilde{\fr{l}} \cap \R Z +  \tilde{\fr{l}}\cap \fr{m}'   
\]
    coincides  with $\fr{l}$. This completes the proof.
  \end{proof}

\subsection{\textsf{Type Ib}}
Now we give a   classification of  minimal   admissible  homogeneous manifolds     $M = G/L$ of \textsf{Type Ib}  with  simple compact $G$. Each such  manifold  has a  reductive decomposition of the  form
\be \label{special contact}  \fr{g} = \fr{l} + \fr{m}=  \fr{l} +(\fr{sp}(1)  + \fr{m}')   \ee
such   that
\[
 \fr{m}^{\fr{l}} = \fr{sp}(1),\quad    C_{\fr{g}}(\fr{sp}(1)) = \fr{l} + \fr{sp}(1) 
 \]
and 
\be \label{special contact1} C_{\fr{g}}(Z) = \fr{l} + \R Z\,,  \quad \text{for  any non-zero} \,\, 
 Z \in \fr{sp}(1)\,.
 \ee
 \noindent It is   easy to  see  that   for  any non-zero $Z \in \fr{sp}(1)$, the  $B$-orthogonal complement $\fr{d} = Z^{\perp} \cap \fr{m} $ defines   an invariant  contact  structure
    on  $M=G/L $.   In particular, the  manifold  $M = G/L$ admits  more than  one  invariant contact structures.
     According to \cite{AS}  such homogeneous manifolds  $M = G/L$ are called  {\it   homogeneous  contact manifolds  of   special  type}.
          It is  proven   that   all   such manifolds are exhausted  by homogeneous spaces of the  form  $M = G/C_G(\fr{sp}(1))$, that is, quotients  of a  compact simple Lie group  modulo the centralizer  of a  regular  3-dimensional  subalgebra $\fr{sp}(1)$. 
   We  get   the  following
   \bp   Minimal admissible homogeneous  manifolds  $M=G/L$ of  a compact simple Lie group $G$   of \textsf{Type  Ib} are   exactly the   contact homogeneous manifolds of special type.
   \ep

    The  following theorem provides a   classification  of  all   such homogeneous spaces. 
    
    \bt\label{WolfFlags}  \cite{AS}  For   any simple  compact Lie group $G$ different from the exceptional Lie group $\G_2$  there  is  a unique  up to isomorphism  minimal admissible  homogeneous  manifold  $$M= G/L = G/C_G(\fr{sp}(1))$$ of \textsf{Type Ib}. It is  the 3-Sasakian manifold   associated   to the  Wolf space   corresponding to $G$, i.e., the homogeneous space
    \[
    W = G/N_G(\fr{sp}(1))= G/\Sp(1) \cdot L\,,
    \]
     where   $\fr{sp}(1) \subset \fr{g}^{\mathbb{C}}$  is the   3-dimensional subalgebra associated to  a  long  root of $\fr{g}^{\bb{C}}$.  When $G=\G_2$, in addition to the just described structure,
        there is   the  quotient  $G/\Sp(1)$, where  $\Sp(1) = C_{\G_2}(\fr{sp}(1))$, and
         $\fr{sp}(1)$ is   the  3-dimensional subalgebra,   associated  to  a short  root.
     \et  
    \begin{proof}
     Let
$$\fr{g}^{\Cc}   = \fr{c} + \sum_{\alpha \in  R} \fr{g}_{\alpha}    $$
be   the   root  space decomposition  of the complexification  $\fr{g}^{\Cc}$ of the  Lie  algebra $\fr{g}$  with respect   to a Cartan subalgebra  $\fr{c}$ of   $\fr{g}^{\Cc}$. We  denote  by $\tau$ the   antiinvolution  (the  complex  conjugation with respect to  $\fr{g}$),   which  determines  the  real form   $\fr{g}$, that is, 
  $\fr{g} = (\fr{g}^{\Cc})^{\tau}$ (see \cite{GOV}).
Then, up to conjugation,  the regular Lie  subalgebra $\fr{s}$ is the   compact  form $\fr{g}(\alpha)^{\tau}$ of  the   subalgebra
\[
\fr{g}(\alpha):= \mathrm{span}(H_{\alpha}, E_{\alpha}, E_{-\alpha})
\]
 associated
to  a  root  $\alpha \in R$ of $\fr{g}^{\bb{C}}$. 
It is well-known that all roots  of the  same length are conjugated.  Therefore, up to conjugation, there  are  either  one   or two   subalgebras $\fr{g}(\alpha)$.  

\noindent Assume first that   $\alpha$ is  a long root   of a  simple  complex Lie algebra  $\fr{g}^{\Cc}$.  Then, the  dual  element  $H_{\alpha }   \in \fr{c}$ defines a (depth 2)  $\Z$-grading
 \[
 \fr{g}^{\Cc} =\fr{g}^{-2} +\fr{g}^{-1}+\fr{g}^{0}+\fr{g}^{1}+\fr{g}^{2} 
 \]
of   $\fr{g}^{\Cc}$,  where  $\fr{g}^j$  is the  eigenspace  of $\ad_{H_{\alpha}}$ with the  eigenvalue $j$. Moreover,
$\fr{g}^{\pm 2} = \Cc E_{\pm \alpha}$,    and  $\fr{g}^{0} = \Cc H_{\alpha} \oplus \fr{l}^{\Cc}$, where  $\fr{l}^{\Cc} = C_{\fr{g}^{\Cc}}({\fr{s}^{\Cc}})$. The   even    subalgebra is  given by
$$\fr{g}^{-2} +  \fr{g}^{0} + \fr{g}^{2}  = \fr{g}(\alpha) \oplus \fr{l}^{\Cc}. $$
It follows that the  condition (\ref{special contact1}) holds true,  and  consequently  $M= G/L$ is  a minimal   homogeneous manifold.

\noindent It  remains  to study the case when $\fr{g}(\alpha)$  is  the   subalgebra associated to a short  root $\alpha$.
    Recall  that  there  are  roots of  different length  only for the simple  Lie groups  of  type ${\sf B}_n, {\sf C}_n,  {\sf F}_4, \G_2$. It is easy to check  that  the root systems $R$ of  ${\sf B}_n, {\sf C}_n$ and  ${\sf F}_4$ satisfy  the  following property:
     for  a short  root $\alpha \in R  $  there  exists  an  orthogonal  to  $\alpha$ root $\beta$ such that   $\alpha + \beta $ is a root. Then,  the root  vector $E_{\beta}$ is  annihilated  by $\ad_{{H}_{\alpha}}$,  but it is not   annihilated  by the  operator $\ad_{E_{\alpha}}$.  In other words,  $E_{\beta} \in C_{\fr{g}^{\Cc}}(H_{\alpha})$, but $E_{\beta } \notin  \fr{l}^{\Cc} = C_{\fr{g}^{\Cc}}(\fr{g}(\alpha))$. This gives a contradiction.

\noindent However, for the   exceptional  Lie algebra $\fr{g}_2$ one can easily check that  the   sum
 $\alpha + \beta$ of  a short  root $\alpha$ and  any   other  root $\beta$ is not a root. This implies  that     the   total spaces   of  both   principal bundles
\[ 
 \pi_1 : \G_2/\Sp(1)^{\alpha}  \to W = \G_2/\Sp(1)^{\alpha}\times \Sp(1)^{\beta}\,,
\]
\[
  \pi_2 : \G_2/\Sp(1)^{\beta}  \to W = \G_2/\Sp(1)^{\alpha}\times \Sp(1)^{\beta}
  \]
  over  the  Wolf  space  $W=\G_2/\Sp(1)^{\alpha}\times \Sp(1)^{\beta}$ are minimally admissible  homogeneous manifolds.  
  \end{proof}
\br
Using the classification of  the Wolf spaces corresponding to  compact simple Lie groups (see, e.g., \cite{Besse}), we list for the convenience of the reader  the admissible decompositions $\fr{g}=\fr{l}+(\fr{sp}(1)+\fr{m}')$ corresponding to the minimal admissible homogeneous spaces associated to the   Theorem \ref{WolfFlags}:
\begin{align*} 
\fr{su}(p+2)=&(\R I+\fr{su}(p))+(\fr{su}(2)+\Cc^p\otimes \Cc^2),\quad I=\diag\left(2iE_p,-piE_2\right)\\
\fr{so}(p+4)=&(\fr{so}(p)+\fr{so}(3))+(\fr{so}(3)+\R^p\otimes \R^4),\\
\fr{sp}(p+1)=&\fr{sp}(p)+(\fr{sp}(1)+\Hh^p),\\
\fr{e}_6=&\fr{su}(6)+(\fr{su}(2)+\wedge^3\Cc^6\,\otimes\Cc^2),\\
\fr{e}_7=&\fr{so}(12)+(\fr{su}(2)+\Delta_{12}\otimes\Cc^2),\\
\fr{e}_8=&\fr{e}_7+(\fr{su}(2)+\wedge^2\R^8\,\otimes\Cc^2),\\
\fr{f}_4=&\fr{sp}(3)+(\fr{su}(2)+\underline{\wedge}^3\Hh^6\,\otimes\Cc^2),\\
\fr{g}_2=&\fr{su}(2)+(\fr{su}(2)+\underline\otimes^3\Cc^2\,\otimes\Cc^2),\\
\fr{g}_2=&\fr{su}(2)+(\fr{su}(2)+\Cc^2\otimes\,\underline\otimes^3\Cc^2),
\end{align*}
where $U \underline\otimes V$ denotes 
denote the highest irreducible component of 
$U\otimes V$, and $\Delta$ denotes the irreducible complex spin module.
\er

 \section{Homogeneous non-compact   Lorentzian manifolds  of  semisimple  Lie groups   with compact stabilizers}\label{secnoncomp}
   Let $G$ be a semisimple  non-compact  Lie group.   We   are interested  in a description of
    homogeneous Lorentzian  manifolds  $M= G/L$   with compact stabilizers $L$. We may  reduce their study  to the description  of all minimal   admissible   homogeneous manifolds $M=G/L$ of this type, or equivalently,   of all maximal  admissible Lie  subgroups $L$ of the group $G$.
      We   consider   only  manifolds  of \textsf{Type  Ia}   and  \textsf{Ib}. Note that there also exist   manifolds  of \textsf{Type Ic} for non-compact semisimple Lie groups.

   \subsection{\textsf{Type Ia}}
    The  study of homogeneous space of this type reduces  to the  description  of all    reductive   decompositions having the  form
    \begin{equation}  
      \label{Type 1a}  \fr{g} = \fr{l }+ \R Z + \fr{m}'\,,
    \end{equation}
    where  $\fr{h} :=  \fr{l} + \R Z = C_{\fr{g}}(Z)$,
    and  $[\fr{h},\fr{m}'] = \fr{m}',$
 and hence to  a description  of  compact  one-dimensional  subalgebras $\R Z$  with compact centralizer    $\fr{h}$.  Let 
 \[
  \fr{g} = \fr{k} + \fr{p} 
  \]
   be  a Cartan decomposition   such that  $\fr{h} \subset \fr{k}$. We   fix   a    Cartan  subalgebra   $\fr{c}_{\fr{k}} $ which contains $Z$   and  extend it to  a Cartan subalgebra  $\fr{c}:= \fr{c}_{\fr{k}}+ \fr{c}_{\fr{p}} $   of $\fr{g}$.
 If    the relation \eqref{Type 1a}   holds,   then   $\fr{c}_{\fr{p}}=0 $  and    any compact  subalgebra  $ \R Z$  with $C_{\fr g} (Z) = \fr{h}$  defines
  a   maximal  admissible    subalgebra $\fr{l}$, which   is the $B$-orthogonal complement  to $Z$  in $\fr{h}$. 
 The condition $\fr{c}_{\fr{p}}=0$ implies 
 $\rk(\fr{k})=\rk(\fr{g})$. First we    consider the  classical   simple real  Lie  algebras  with  a compact Cartan subalgebra $\fr{c} \subset \fr{k}$, with aim to  describe   all  elements   $Z \in \fr{c}$
 having  trivial centralizers    $C_{\fr{p}} (Z)=0$.  Then,  $\fr{l}$ will be the orthogonal complement of $Z$ in $C_{\fr{k}} (Z)=0$, $\fr{m}'$ will be the orthogonal complement to $\fr{l}$ in $\fr{g}$, and moreover $\fr{m}=\R Z+\fr{m}'$ (direct sum). We use these facts, to obtain a description for 
   all classical   semisimple real  Lie  algebras.
   
   \br
First it is convenient to pose   the Cartan decompositions $ \fr{g} = \fr{k} + \fr{p} $
of all classical real simple non-compact Lie algebras $ \fr{g}$ with $\rk\fr{k}=\rk\fr{g}$ (we follow   \cite{GOV}).
\begin{align*}
\fr{su}(p,q)&=s\big(\fr{u}(p)+\fr{u}(q)\big)+\Cc^{p}\otimes\Cc^{q},\quad p, q\geq 1,\\
\fr{so}(p,q)&=\big(\fr{so}(p)+\fr{so}(q)\big)+\R^{p}\otimes\R^{q},\quad p, q\geq 1,\quad \text{$pq$ is even},\\
\fr{sp}(2n,\R)&=\fr{u}(n)+S^2\Cc^n,\quad n\geq 1,\\
\fr{sp}(p,q)&=\big(\fr{sp}(p)+\fr{so}(q)\big)+\Hh^{p}\otimes_{\Hh}\Hh^{q},\quad p, q\geq 1,\\
\fr{so}(n,\Hh)&=\fr{u}(n)+\wedge^2\Cc^n,\quad n\geq 1. 
\end{align*}
Now we are ready to consider these algebras case-by-case.
\er

$\bullet$ Let $\fr{g}=\fr{su}(p,q)$, $ p, q\geq 1$. Then $$\fr{g}=\fr{k}+\fr{p}=s\big(\fr{u}(p)+\fr{u}(q)\big)+\Cc^{p}\otimes\Cc^{q},$$ where for a subalgebra $\fr{f}\subset\fr{u}(m)$, $s(\fr{f})$ denotes the intersection $\fr{f}\cap\fr{su}(m)$.
The Lie algebra $\fr{su}(p,q)$ consists of the complex matrices
$$\left(\begin{matrix} X_1& Y\\ \bar Y^T&  X_2\end{matrix}\right),\quad \bar X_1^T=-X_1,\quad \bar X_2^T=-X_2,\quad \tr X_1+\tr X_2=0,$$
here $X_1$, $X_2$, and $Y$ are blocks of the size $p\times p$, $q\times q$, and $p\times q$, respectively. 
Elements of $\fr{k}$ are given by the matrices with $Y=0$, and elements of $\fr{p}$ are given by the matrices with $X_1=X_2=0$.
The standard Cartan subalgebra $\fr{c}=\fr{c}_\fr{p}$ consists of  elements of $\fr{k}$ with diagonal matrices $X_1,X_2$. 
Fix some $Z\in\fr{c}$ and let 
\[
\Cc^p=\oplus_{i=1}^r V_i,\quad  \Cc^q=\oplus_{\alpha=1}^s U_\alpha
\]
be the eigenspace decompositions of the endomorphisms $Z|_{\Cc^p}$ and  
$Z|_{\Cc^q}$, respectively.
Then 
\[
Z=\sum_{j=1}^rb_jI_j+\sum_{\alpha=1}^s c_\alpha I'_\alpha\,,
\]
where $b_j$, $c_\alpha\in\R$, and $I_j$ is an endomorphism acting on $V_j$ as the multiplication by $i$ and annihilating the orthogonal complement to $V_j$; the endomorphisms $I'_\alpha$ are defined in a similar way. Since $Z\in\fr{su(p+q)}$, it holds 
\[
\sum_{j=1}^r b_j \dim_{\Cc} V_j+\sum_{\alpha=1}^s c_\alpha\dim_{\Cc}  V_\alpha=0\,.
\]
  Since $\exp(tX)$ is compact, the numbers $b_{1},\dots, b_{p},c_{1},\dots, c_{q}$ are commensurable, and we may assume that these numbers are rational. The condition $C_{\fr{m}'}(Z)=0$ implies $C_{\fr{p}}(Z)=0$. Let an element $A\in \fr{p}$ be given by the matrix $Y=E_{i\alpha}$, where $E_{i\alpha}$ is the matrix with 1 at the position $(i,\alpha)$ and zeros in the rest  entries. Then $[Z,A]$ is given by the matrix $(b_j-c_\alpha)E_{j\alpha}$. This implies that the condition $C_{\fr{p}}(Z)=0$ is equivalent to the condition
$$b_{j}\neq c_{\alpha},\quad 1\leq j\leq r,\quad 1\leq \alpha\leq s.$$
Now we see that
$$C_\fr{k}(Z)=s\big(\oplus_{i=1}^r \fr{u}(V_i)\,+\,  \oplus_{\alpha=1}^s \fr{u}(U_\alpha)  \big)$$
and  $$\fr{l}=\oplus_{i=1}^r \fr{su}(V_i)\,+\,  \oplus_{\alpha=1}^s \fr{su}(U_\alpha)+\hat {\fr{l}},$$ where $\hat{\fr{l}}$
is the orthogonal complement to $Z$ in $s(\left<I_1,\dots,I_r,I'_1,\dots,I'_s \right>)$ with repsect to the Killing form.
Finally we obtain the decomposition
$$\fr{g}=\fr{l}+\fr{m}=\fr{l}+(\R Z+\fr{m}'),$$
where $$\fr{m}'=\fr{p}+\sum_{1\leq j<k\leq r} V_j\otimes_{\Cc} V_k+\sum_{1\leq\alpha<\beta\leq s} U_\alpha\otimes_{\Cc} U_\beta.$$

$\bullet$  Let $\fr{g}=\fr{so}(p,q)$, $ p, q\geq 1$. Then $$\fr{g}=\fr{k}+\fr{p}=\big(\fr{so}(p)+\fr{so}(q)\big)+\R^{p}\otimes\R^{q}.$$
The condition $\rk(\fr{g})=\rk(\fr{k})$ is fulfilled whenever at least on of the numbers $p$, $q$ is even.  Fix orthonormal bases $e_1,\dots,e_p$ and  $f_1,\dots,f_q$ of the spaces $\R^p$ and $\R^q$, respectively.
The Cartan subalgebra $\fr{c}$ of $\fr{so}(p,q)$ is the direct sum of the Cartan subalgebras of $\fr{so}(p)$ and $\fr{so}(q)$. The Cartan subalgebra of $\fr{so}(p)$ consists of the elements 
$$e_1\wedge e_2,\dots, e_{[\frac{p}{2}]-1}\wedge e_{[\frac{p}{2}]}.$$
A similar structure has the Cartan subalgebra of  $\fr{so}(q)$.
For a given element $Z\in\fr{c}$ we may assume that the basis is adapted to the canonical form of the element $Z$, i.e., there are decompositions
 $$\R^p=\oplus_{i=1}^r V_0,\quad  \R^q=\oplus_{\alpha=0}^s U_\alpha$$
such that $$Z=\sum_{j=1}^rb_jI_j+\sum_{\alpha=1}^s c_\alpha I'_\alpha,$$
where $b_j$, $c_\alpha\in\R$, and $I_j$, $I'_\alpha$ are  complex structures on $V_j$, $U_\alpha$, respectively. As in the previous case, we may assume that 
the numbers $b_1,\dots,b_r,c_1,\dots,c_s$ are rational; the condition $\fr{c}_{\fr{p}}(Z)=0$ implies 
$$b_{j}\neq c_{\alpha},\quad 1\leq j\leq r,\quad 1\leq \alpha\leq s,$$
and one of the spaces $V_0$ or $U_0$ is trivial. We assume that $U_0=0$.
 Moreover, if $q$ is odd, then $V_0=0$.
It is clear that
$$C_\fr{k}(Z)=\fr{so}(V_0)+\oplus_{i=1}^r \fr{u}(V_i)\,+\,  \oplus_{\alpha=1}^s \fr{u}(U_\alpha),$$ and
 $$\fr{l}=\fr{so}(V_0)+\oplus_{i=1}^r \fr{su}(V_i)\,+\,  \oplus_{\alpha=1}^s \fr{su}(U_\alpha)+\hat {\fr{l}},$$ where
 $\hat {\fr{l}}$ is the orthogonal complement to $Z$  in $\left<I_1,\dots,I_r,I'_1,\dots,I'_s \right>$ with respect to the Killing form.
 Recall that there is a 2-grading $$\fr{so}(2m)=\fr{u}(m)+\wedge^2\Cc^{m}.$$
 This implies that
 $$\fr{g}=\fr{l}+\fr{m}=\fr{l}+(\R Z+\fr{m}'),$$
where $$\fr{m}'=\fr{p}+\sum_{j=1}^r\wedge^2_{\Cc}V_j+\sum_{\alpha=1}^s\wedge^2_{\Cc}U_\alpha+\sum_{0\leq j<k\leq r} V_j\otimes V_k+\sum_{1\leq\alpha<\beta\leq s} U_\alpha\otimes U_\beta.$$

$\bullet$ Let $\fr{g}=\fr{sp}(2n,\R),$ $n\geq 1$. Then 
$$\fr{g}=\fr{k}+\fr{p}=\fr{u}(n)+S^2\Cc^n.$$
The Lie algebra $\fr{sp}(2n,\R)$ consists of the real matrices
$$\left(\begin{matrix} X& Y_1\\  Y_2&  -X^T\end{matrix}\right),\quad \bar Y_1^T=Y_1,\quad \bar Y_2^T=Y_2,$$
here $X$, $Y_1$, and $Y_2$ are square $n\times n$ matrices. 
Elements of $\fr{k}$ are given by the matrices with $X^T=-X$ and $Y_2=-Y_1$, while elements of $\fr{p}$ are given by the matrices with $X^T=X$ and $Y_2=Y_1$.
The standard Cartan subalgebra $\fr{c}=\fr{c}_\fr{p}$ consists of the elements of $\fr{k}$ given by diagonal matrices $Y_1$. 
Let $Z\in\fr{c}$ be given by the matrix $Y_1=\diag(z_1,\dots, z_n)$. The condition $C_{\fr{p}}(Z)=0$  implies easily
that the numbers $z_i$ are non-zero and pairwise different.
The center $C_{\fr{k}}(Z)$ coincides with $\fr{c}$. Then $\fr{l}$ is made up of the elements from $\fr{c}$ given by $Y=\diag(y_1,\dots,y_n)$ satisfying the condition $y_1z_1+\cdots+ y_nz_n=0$. Finally,
$$\fr{g}=\fr{l}+\fr{m}=\fr{l}+(\R Z+\fr{m}'),$$
where $$\fr{m}'=\fr{p}+\left\{\left(\begin{matrix} X& Y\\  -Y&  X\end{matrix}\right) \in\fr{u}(n)| Y \text{ has zeros on the diagonal}  \right\}.$$

$\bullet$ Let $\fr{g}=\fr{sp}(p,q)$, $ p, q\geq 1$. Then $$\fr{g}=\fr{k}+\fr{p}=\big(\fr{sp}(p)+\fr{so}(q)\big)+\Hh^{p}\otimes_{\Hh}\Hh^{q}.$$ 
The elements of the Lie algebra $\fr{sp}(p,q)$ may be identified with the quaternionic matrices $$\left(\begin{matrix} X_1& Y\\ \bar Y^T&  X_2\end{matrix}\right),\quad \bar X_1^T=-X_1,\quad \bar X_2^T=-X_2,$$
here $X_1$, $X_2$, and $Y$ are blocks of the size $p\times p$, $q\times q$, and $p\times q$, respectively. 
Elements of $\fr{k}$ are given by the matrices with $Y=0$, and elements of $\fr{p}$ are given by the matrices with $X_1=X_2=0$.
The standard Cartan subalgebra $\fr{c}=\fr{c}_\fr{p}$ consists of the elements of $\fr{k}$ with diagonal matrices $X_1,X_2$ with imaginary complex numbers on the diagonal. 
Let $Z\in\fr{c}$. It is given by the diagonal matrix
$$\diag(0E_{p_0},z_1iE_{p_1},\dots, z_ri E_{p_r},0E_{q_0},z'_1iE_{q_1},\dots, z'_s iE_{q_s}),$$
where $z_i$, $z'_\alpha$ are real numbers; the numbers $z_1,\dots,z_r$ are pairwise different; the same holds for the numbers $z'_1,\dots,z'_s$. The condition  
 $C_{\fr{p}}(Z)=0$ implies
$$z_{j}\neq z'_{\alpha},\quad 1\leq j\leq r,\quad 1\leq \alpha\leq s,$$
and at least one of the numbers $p_0$, $q_0$ is zero. We assume that $q_0=0$.
Recall that there is a decomposition
$$\Hh^m=\Cc^m+j \Cc^m,$$
and it holds $$C_{\fr{sp}(m)}(\operatorname{Op}(iE_m))=\fr{u}(m),$$
where  $\operatorname{Op}(iE_m)$ is an element of $\fr{sp}(m)$ with the matrix $iE_m$, and $\fr{u}(m)$ acts diagonally in $\Hh^m=\Cc^m+j \Cc^m$. 
Moreover, there is a $\Z$-grading
$$\fr{sp}(m)=\fr{u}(m)+S^2\Cc^m.$$

This implies that
$$C_\fr{k}(Z)=\fr{sp}(p_0)\,+\,\oplus_{i=1}^r \fr{u}(p_i)\,+\,  \oplus_{\alpha=1}^s \fr{u}(q_i)  $$
and  $$\fr{l}=\fr{sp}(p_0)\,+\, \oplus_{i=1}^r \fr{su}(p_i)\,+\,  \oplus_{\alpha=1}^s \fr{su}(q_i)  +\hat {\fr{l}},$$ where $\hat{\fr{l}}$
is the orthogonal complement to $Z$ in the space spanned by operators corresponding to the matrices $iE_{p_1},\dots, i E_{p_r},iE_{q_1},\dots, iE_{q_s}$.
We obtain the decomposition
$$\fr{g}=\fr{l}+\fr{m}=\fr{l}+(\R Z+\fr{m}'),$$
where $$\fr{m}'=\fr{p}+\sum_{j=1}^rS^2\Cc^{p_j}+\sum_{\alpha=1}^sS^2\Cc^{q_\alpha}
+ \sum_{0\leq j<k\leq r} \Hh^{p_j}\otimes_{\Hh}{\Hh}^{p_k} +\sum_{1\leq\alpha<\beta\leq s} \Hh^{q_\alpha}\otimes_{\Hh} \Hh^{q_\beta}.$$

 $\bullet$ Let $\fr{g}=\fr{so}(n,\Hh),$ $n\geq 1$. Then 
 $$\fr{g}=\fr{k}+\fr{p}=\fr{u}(n)+\wedge^2\Cc^n.$$
 The Lie algebra $\fr{so}(n,\Hh)$ consists of the complex matrices
 $$\left(\begin{matrix} X& Y\\  -\bar Y&  \bar X\end{matrix}\right),\quad \bar X^T=-X,\quad \bar Y^T= Y,$$
 here $X$, $Y$ are square $n\times n$ matrices. 
 Elements of $\fr{k}$ are given by the matrices with $\bar X=X$ and $\bar Y=Y$, while elements of $\fr{p}$ are given by the matrices with $\bar X=-X$ and $\bar Y=-Y$.
 The standard Cartan subalgebra $\fr{c}=\fr{c}_\fr{p}$ consists of the elements of $\fr{k}$ given by diagonal matrices $X$ with imaginary elements at the diagonal. 
 Let $Z\in\fr{c}$ be given by the matrix $$X=\diag(z_1i,\dots, z_ni).$$ The condition $C_{\fr{p}}(Z)=0$ easily implies
 that the numbers $z_i$ are non-zero and pairwise different.
 The center $C_{\fr{k}}(Z)$ coincides with $\fr{c}$. Then $\fr{l}$ is made up of the elements from $\fr{c}$ given by $V=\diag(v_1i,\dots,v_ni)$ satisfying the condition $v_1z_1+\cdots+ v_nz_n=0$.
  Finally,
 $$\fr{g}=\fr{l}+\fr{m}=\fr{l}+(\R Z+\fr{m}'),$$
 where $$\fr{m}'=\fr{p}+\left\{\left(\begin{matrix} X& Y\\  -Y&  X\end{matrix}\right) \in\fr{u}(n)| Y \text{ has zeros at the diagonal}  \right\}.$$
 
 \medskip
 
 \medskip
 
 Suppose now that $\fr{g}=\fr{g_1}\oplus\cdots\oplus\fr{g}_r$ is a real semisimple Lie algebra such that its simple ideals $\fr{g}_i$  are classical non-compact Lie algebras.  
  As before consider the decomposition   $$\fr{g}=\fr{l}+\fr{m}=\fr{l}+(\R Z+\fr{m}').$$
  The vector $Z$ may be represented as
  $$Z=Z_1+\cdots+Z_r, \quad Z_i\in\fr{g}_i.$$
The maximality condition implies that  all $Z_i$ are non-zero.
We immediately conclude that 
$$\fr{g}=\fr{l}+\fr{m}=\left( \sum_{i=1}^r\fr{l}_i+\hat{\fr{l}}\right)+\left(\R Z+\sum_{i=1}^r\fr{m}'_i\right),$$
where $\fr{g}_i=\fr{l}+(\R Z_i+\fr{m}_i')$ are the just obtained decompositions, and $\hat{\fr{l}}$ is the orthogonal complement to $Z$ in $\left<Z_1,\dots,Z_r\right>$.
 
\subsection{\textsf{Type Ib}}
In this case the problem    reduces  to the description of homogeneous manifolds $M = G/L$  of   non-compact simple Lie groups $G$  with compact stabilizers $L$,    which admit  a reductive decomposition   of the  form
\[
    \fr{g} = \fr{l} + \fr{m} =  \fr{l}  + (\fr{s} +\fr{m}')
    \]
where $ \fr{h} = \fr{l}+\fr{s}$ (direct sum)   is the  normalizer $N_{\fr{g}}(\fr{s})$
of   a  3-dimensional simple   regular   subalgebra isomorphic to  $\fr{sp}(1)$ or  $\fr{sl}(2,\R)$.  The   complexification
$\fr{g}^{\Cc}  = \fr{l}^{\Cc} + \fr{s}^{\Cc} +(\fr{m}')^{\Cc}    $
is the  standard  decomposition   associated to the  3-dimensional regular subalgebra
$$\fr{s}^{\Cc} = \fr{g}(\alpha) = \mathrm{span} (H_{\alpha}, E_{\pm \alpha})$$
of   $\fr{g}^{\Cc}$,  associated  to  a  long  root  $\alpha  \in R$, and  for the case of the Lie algebra  of $\fr{g}_2$   also to a short   root, see  Section \ref{compactcase}.
   The complex subalgebra $\fr{g}(\alpha)$  generates  a  subgroup $G(\alpha)$
    of the Lie group   $G^{\Cc}$,   and   the  homogeneous manifold 
    \[
    M^{\Cc} = G^{\Cc}/ N_{G^{\Cc}}(G(\alpha)) =  G^{\Cc}/ H^{\Cc}
    \]
     is the  complexification
     of the   Wolf  space  $   (G^{\Cc})^{\tau}/ (H^{\Cc})^{\tau}$ corresponding to the   compact
      form  $(G^{\Cc})^{\tau}$  of  $G^{\Cc}$.
   The    classification of  such decompositions for   real simple Lie  algebras $\fr{g} $   reduces 
    to  the classification  of all  real   forms $\fr{g} $  of  the  complex   simple Lie  algebras   $\fr{g}^{\Cc}$, defined  by an  anti-involution
   $\sigma $  which  preserves  the  subalgebra $\fr{g}(\alpha)$.   Then, the  fixed point    set  $(\fr{g}(\alpha))^{\sigma}$  of the    subalgebra   $\fr{g}(\alpha)$ is either   the compact real  form
    $\fr{sp}(1)$  or the   non-compact real form   $\fr{sl}(2,\R)$.
    Hence, our task  is   equivalent to  the  classification   of   pseudo-Riemannian symmetric quaternionic-K\"ahler  manifolds   $Q= G/\Sp(1) \cdot L$   solved  in \cite{AleC05},  and para-quaternoinic K\"ahler symmetric   manifolds  $P = G/\Sl(2, \R)\cdot L$, solved  in \cite{C}  and \cite{DJS}.
    
      Select a non-compact   homogeneous manifold $M = G/L$     with   compact  stabilizer $L = C_{G}(A_1)$, where  $A_1$ is isomorphic to  $\Sp(1)$ or to $\Sl(2,\R)$.
    Assume that   $\fr{g}(\alpha)$  is the  3-dimensional  Lie algebra  associated to a long root. In  the  first  case  the homogeneous  manifold $Q = G/H$,  where  $G$  is the Lie group  with  the Lie   algebra $\fr{g} = (\fr{g}^{\Cc})^{\tau}$
and $H = L \cdot\Sp(1)$  is  the  subgroup  generated  by  $\fr{h} = (\fr{l}^{\Cc})^{\sigma} + \fr{g}(\alpha)^{\tau}$, is a  symmetric  pseudo-quaternionic K\"ahler  manifold. All such  manifolds  are  classified in  \cite{AleC05}. 
In the second case, the manifold $Q= G/H = G/L \cdot\Sl(2, \R)  $ is a  symmetric para-quaternionic  K\"ahler  manifold. All such manifolds  are classified in \cite{DJS}. 
Now, we need the  classification of  all
anti-involutions $\sigma$  of the  complex  simple  Lie algebra
 $\fr{g}^{\Cc}$,  with the  standard  grading
 \[
   \fr{g}^{\Cc} = \fr{g}^{-2} + \fr{g}^{-1} + \fr{g}^{0}+ \fr{g}^1 + \fr{g}^2
   \]
 preserved by   the  grading   element $H_{\mu} \in \fr{c}$. Such  real forms   had been  classified   in \cite{C}.
 To get our   classification, we  have   to examine  the  list  of all pseudo-Riemannian  symmetric  quaternionic K\"ahler spaces $G/\Sp(1) \cdot L$ presented in  \cite{AS}  and the list of  all para-quaternionic-K\"ahler  symmetric  spaces given in
 \cite{DJS}  or \cite{C}, and from the stated homogeneous spaces     select the  manifolds  with  compact stabilizer $L$.

  In the first  case, we deduce that $M=G/L $ is  the  total space of the  principal $\Sp(1)$-bundle
\[
 \pi :M= G/L \to Q=  G/L \cdot\Sp(1)
 \]
  which is a pseudo-Riemannian 3-Sasakian manifold.   In the second case, $M$ is the   total space   of the  principal $\Sl(2,\R)$-bundle
\[
 \pi :M= G/L \to Q=  G/L \cdot\Sl(2,\R),
 \]
 which is a  para-3-Sasakian manifold.
  
  \br
   Note that   the    complex  exceptional  Lie    group
$\G_2^{\Cc}$  has  two  real forms. The  normal  form $G_{2(2)}$ with  maximal compact  subgroup  $K = \Sp(1)^{sh} \cdot\Sp(1)^{l}$  corresponding  to  short  root   and  orthogonal to  it long root,   and  the compact form  $\G_2$,  which  defines  the     quaternionic-K\"ahler   symmetric   space     $\G_2/\Sp(1)^{sh} \cdot\Sp(1)^{l}$.  Hence,  there   are   two minimal
admissible    manifolds of \textsf{Type Ib}  of the non-compact  group $\G_{2(2)}$, namely:   $\G_{2(2)}/\Sp(1)^{sh}$  and
$\G_{2(2)}/\Sp(1)^{l}$.
\er

We may summarize as follows:

\bt \label{Para} Let $G$ be a non-compact simple Lie group.  All  minimal  admissible non-compact   homogeneous manifolds $M = G/L$   of    \textsf{Type Ib} corresponding to $G$ are exhausted  by   pseudo-3-Sasakian  manifolds $M = G/L = G/C_{G}(\fr{sp}(1))$, by   para-3-Sasakian manifolds  $M = G/L = G/C_{G}(\fr{sl}(2, \R))$, both associated to a long root, and  also  by the non-compact   homogeneous  manifold $\G_2/\Sl(2,\R)$, where $\Sl(2, \R)$ is the subgroup  corresponding to a  simple Lie  algebra  associated  to a short root of the  split  Lie algebra $\fr{g}_{2(2)}$.
  \et
  
\br
 Using the classification given in \cite{AleC05}, we may present  the admissible decompositions $\fr{g}=\fr{l}+(\fr{sp}(1)+\fr{m}')$ corresponding to the homogeneous spaces of Theorem \ref{Para}:
 \begin{align*} 
 \fr{su}(p,2)=&(\R I+\fr{su}(p))+(\fr{su}(2)+\Cc^p\otimes \Cc^2),\quad I=\diag\left(2iE_p,-piE_2\right)\\
 \fr{so}(p,4)=&(\fr{so}(p)+\fr{so}(3))+(\fr{so}(3)+\R^p\otimes \R^4),\\
  \fr{sp}(p,1)=&\fr{sp}(p)+(\fr{sp}(1)+\Hh^{p,1}),\\
 \fr{e}_{6(2)}=&\fr{su}(6)+(\fr{su}(2)+\wedge^3\Cc^6\,\otimes\Cc^2),\\
 \fr{e}_{7(-5)}=&\fr{so}(12)+(\fr{su}(2)+\Delta_{12}\otimes\Cc^2),\\
 \fr{e}_{8(-24)}=&\fr{e}_{7}+(\fr{su}(2)+\wedge^2\R^8\,\otimes\Cc^2),\\
 \fr{f}_{4(4)}=&\fr{sp}(3)+(\fr{su}(2)+\underline{\wedge}^3\Hh^6\,\otimes\Cc^2),\\
 \fr{g}_{2(2)}=&\fr{su}(2)+(\fr{su}(2)+\underline\otimes^3\Cc^2\,\otimes\Cc^2),\\
 \fr{g}_{2(2)}=&\fr{su}(2)+(\fr{su}(2)+\Cc^2\otimes\,\underline\otimes^3\Cc^2).
 \end{align*}
 Based on the results  in \cite{DJS}, we deduce  that  there is only one admissible decomposition\\ $\fr{g}=\fr{l}+(\fr{sl}(2,\R)+\fr{m}')$:
 $$\fr{su}(p+1,1)=(\R I+\fr{su}(p))+(\fr{su}(1,1)+\Cc^p\otimes \Cc^{1,1}),\quad I=\diag\left(2iE_p,-piE_2\right).$$
\er

\section{Homogeneous Lorentzian manifolds with a reductive stabilizer   of  \textsf{Type  II} and \textsf{Type III}}\label{IIIII}
In this final section we present the classification of homogeneous Lorentzian manifolds $G/L$ with a reductive stabilizer   of  \textsf{Type  II} and \textsf{Type III} (under certain assumptions).  Let us first  treat the case corresponding to \textsf{Type  II} stabilizers.
\bt \label{typeb}
Let $M =  G/L$  be  a simply   connected  (almost)  effective  homogeneous  Lorentz manifold  with  isotropy  subalgebra of  $\mf{Type \  II}$,  that  is  $j_{*}(\fr{l}) = \fr{so}(W)+ \fr{k}$,  where  $W$ is  a Lorentzian subspace of  dimension  $m =\dim W >2$ and $\fr{k}$ is  a  compact  subalgebra  of  $\fr{so}(W^{\perp})$. Suppose that  $\fr{k}$ does not annihilate  any non-zero vector in $W^{\perp}$.  Then, 
  $M$    must be  a  direct product  of a $m$-dimensional homogeneous Lorentz  space $M_0$ and of a   homogeneous  Riemannian manifold  $\mc{N}= G_1/N$. In particular,
  \begin{itemize}
  \item if $m>3$, then  $M_0$ is a space of  constant curvature, that is    the Minkowski  space $\mf{M}^{1, m-1}=\R^{1,m-1}$, or  the de Sitter  space  $\mf{dS}^{m}= \SO(1,m)/\SO(1,m-1)$,  or  the
  anti de Sitter  space $\mf{AdS}^{m}= \SO(2,m-1)/\SO(1,m-1)$.
\item   if $m=3$, then $M_0$ is either the space of constant curvature, or the Lie group $\widetilde{\rm{SL}(2,\mathbb{R})}$.
\end{itemize}
\et

\begin{proof} When $L$ is of \textsf{Type  II}, then the   reductive  decomposition  can be  written  as
	\[
	\fr{g} = \fr{l} + \fr{m}  = (\fr{so}(W)+ \fr{k}) + (W +\mc{U})
	\]
	where  $\mc{U} = W^{\perp}$ is  an Euclidean  vector  space   and  $\fr{k} \subset \fr{so}(\mc{U})$.
	The map
	\[
	\Lambda^2 (W + \mc{U})=  \Lambda^2 W +  W \otimes \mc{U} + \Lambda^2\mc{U}\to  \fr{g}\,,
	\]
	given by the Lie brackets,  is  a $(\fr{so}(W)+\fr{k})$-equivariant linear map. This implies the relations
	\[
	[W,W]\subset\fr {so}(W)+W\,,\quad
	   [W, \mc{U}] =0\,,\quad  [\fr{k},\mc{U}]\subset \mc{U}\,.
	   \]
	    Thus, $\fr{g}$ is the direct sum of two ideals, i.e.,
	    \[
	    \fr{g} =  (\fr{so}(W)+ W) + (\fr{k} +\mc{U})\,.
	    \]
	Hence, $M$ must be a product of a homogeneous Lorentzian space $M_0$ and a homogeneous Riemannian manifold. This proves our first claim. \\
	Consider now  the Lie algebra $\fr{so}(W)+ W$. If $\dim W>3$, then
		the  Lie  bracket  restricted to $W$  is  given  by
	$ [w_1, w_2 ] = cw_1 \wedge w_2 \in \fr{so}(W)$ for some $c\in\mathbb{R}$. Suppose that $\dim W=3$. Then, there is an isomorphism $\varphi: \fr{so}(1,2)\to\mathbb{R}^{1,2}$ of the $\fr{so}(1,2)$-modules. For the Lie brackets of $X,Y\in\mathbb{R}$ it holds
	\[
	[X,Y]=c_1 X\wedge Y+c_2 \varphi (X\wedge Y)\,.
	\]
	If $c_1\neq 0$, then the reductive decomposition may be chosen in such a way that $c_2=0$, and  then $\fr{so}(W)+ W$ determines a symmetric space. Suppose that $c_1=0$ and $c_2\neq 0$. Then $\fr{so}(W)+ W$ determines a flat naturally reductive homogeneous space $M_0$. Since $[W,W]\subset W$ is given by the Lie brackets on $\fr{so}(1,2)$ and $M_0$ is simply connected, we conclude that $M_0$ is the Lie group
	$\widetilde{\rm{SL}(2,\mathbb{R})}$.
	   This completes the proof.
	\end{proof}
	
	   \bt 
	   Let $M =  G/L$  be  a simply   connected  (almost)  effective  homogeneous  Lorentz manifold  with  isotropy  subalgebra of  $\mf{Type \  III}$,  that  is  $j_{*}(\fr{l}) = \mathbb{R} d+ \fr{k}$,   $d=p\wedge q+C_0$, $C_0\in\fr{so}(E)$, and $\fr{k}$ is  a  compact  subalgebra  of  $\fr{so}(E)$; $\fr{k}$ commutes with $C_0$. Suppose that  $\fr{k}$ does not annihilate  any non-zero vector in $E$.  Then
	   $M$    is  a  direct product  of    $2$-dimensional constant curvature Lorentz  space $M_0$	
	   and of a   homogeneous  Riemannian manifold  $\mc{N}= G_1/N$.
	   \et

	\begin{proof}
		Since the restriction of the Lie bracket $\Lambda^2 \fr{m}\to\fr{g}$ is an $(\mathbb{R} d+\fr{k})$-equivariant map, we conclude that
		\[
		[p,q]=\lambda d\,, \   (\lambda\in\mathbb{R})\,,\quad [p,E]=[q,E]=0\,,\quad [E,E]\subset \mathbb{R}d+\fr{k}+E\,.
		\]
		Then, the Jacobi identity yields the relation $[[X,Y],p]=0$, for any $X,Y\in E$, i.e.,
		$[E,E]\subset \fr{k}+E$. Similarly, $[[p,q],X]=0$, i.e.,
		either $\lambda=0$ or $C_0=0$. This easily implies the proof.
		\end{proof}
Let us finally highlight the following remark about the subspace $E_0\subset E$ consisting of vectors annihilated by $\fr{k}$, and its triviality when one uses maximal admissible subgroups.
\br Let $M =  G/L$  be  a simply   connected  (almost)  effective  homogeneous  Lorentz manifold  with  isotropy  subalgebra of  \textsf{Type II} or \textsf{Type III}. Suppose that the Lie subalgebra $\fr{l}\subset\fr{g}$ is maximal admissible.  In this case we claim that the subspace $E_0\subset E$ consisting of vectors annihilated by $\fr{k}$ is trivial. Indeed, it is easy to see that the decomposition
\[
\fr{g}=(\fr{l}+E_0)+(\mathbb{R}p+E_1+\mathbb{R}q)
\]
  contradicts the maximality assumption (here $E_1$ is the orthogonal complement to $E_0$ in $E$).
\er

\end{document}